# Detection of coarse-grained unstable states of microscopic/stochastic systems: a timestepper-based iterative protocol

A. C. Tsoumanis and C. I. Siettos[*]

*School of Applied Mathematics & Physical Sciences,
National Technical University of Athens, GR 157 80, Athens, Greece*

**Abstract.** We address an iterative procedure that can be used to detect coarse-grained hyperbolic unstable equilibria (saddle points) of microscopic simulators when *no equations at the macroscopic level are available*. The scheme is based on the concept of *coarse timestepping* [Kevrekidis et al., 2003] and treats the detailed available simulator as an evolving experiment. The iterative protocol incorporates an adaptive mechanism based on the chord method allowing the location of coarse-grained saddle points directly. Ultimately, it can be used in a consecutive manner to trace the coarse-grained open-loop saddle-node bifurcation diagrams. We illustrate the procedure through two indicatively examples including (i) a kinetic Monte Carlo simulation (kMC) of simple surface catalytic reactions and (ii) a simple agent-based model, a financial caricature which is used to simulate the dynamics of buying and selling of a large population of interacting individuals in the presence of mimesis. Both models exhibit coarse-grained regular turning points which give rise to branches of saddle points.

## 1. Introduction

Many problems of current engineering, physical and biological interest are characterized, due to the complexity of the underlying physics, by the lack of good macroscopic models in the form of ordinary, partial and/or integro-differential equations that can be used to describe the dynamic behaviour in the continuum. When such continuum mathematical models are not available, microscopic simulations constitute the backbone of complex systems analysis practice. Several techniques ranging from kinetic Monte Carlo (kMC) to Molecular and Brownian Dynamics and from Cellular Automata to Individual and Agent-based simulations are being used extensively in many research areas as diverse as materials science [1, 2, 3, 4, 5, 6, 7, 8, 9], fluid mechanics [10, 11, 12, 13, 14, 15, 16], socio-economic fields [17, 18, 19, 20, 21], neuroscience [22, 23, 24, 25, 26, 27, 28, 29], epidemiology [30, 31, 32, 33, 34], ecology [35, 36, 37], biology [38, 39], car traffic and pedestrian dynamics management [40, 41, 42, 43, 44].

However, due to the inherent stochasticity and the many-particles (atoms, molecules, macro-molecules, cells, tissues, individuals, populations) interactions, the emergent systems behavior at the macroscopic-continuum level is far from trivial to predict. Under this perspective, detailed simulations are used as "computer experiments" in (a "loose") analogy to the "physical experiments": one sets up the (microscopic/experimental) initial conditions, the values of the parameters and performs runs (leaves the experiment to evolve) for long times (and for many microscopic ensembles consistent to the same initial conditions) in order to observe the macroscopic behaviour. Yet, this is but the first thing one may employ in order to analyse the system's coarse-grained dynamics. Important tasks such as the exact location of the critical points that mark the onset of instabilities, cannot be easily obtained just through temporal simulations.

For the systematic and effective analysis of the simulation dynamics one has to resort to bifurcation analysis. Numerical bifurcation theory provides an arsenal of algorithms and software packages [45, 46, 47, 48, 49, 50]. While these are invaluable tools for performing systematic analysis for small to medium scale systems, there are some drawbacks in using them. On one hand, most of them require from the user to input the system evolution equations, which are assumed to be explicitly available in a specific way. One the other hand, linking them with legacy/commercial or "home-made" developed simulators is not something trivial. Furthermore, the convergence on unstable equilibria through such Newton-like solvers is most of the times impossible when comes to microscopic/stochastic systems. In particular, for microscopic models, the lack of accurate macroscopic representations constitutes a major obstacle towards the systematic investigation of the systems behaviour at the macroscopic level where the analysis is usually sought for design purposes.

Over the last years, it has been established that the so called "Equation-Free" approach may efficiently confront with the problem in an integrated way [51, 52, 53, 54, 55, 56, 57, 58, 59, 60]. This multi-scale model reduction framework can be exploited to perform coarse-grained tasks such as bifurcation and stability analysis [51, 53, 54, 55, 61, 57, 58, 59], controller design [62, 63], computation of travelling waves and periodic solutions [64, 65], in a wide range of complex problems in mechanics and engineering, bypassing the need for extracting continuum-level equations in a closed form. The available detailed microscopic simulator is treated as a black-box timestepper of the coarse-grained observables; vital coarse-grained quantities, required for continuum numerical analysis (coarse residuals, the action of coarse Jacobians, control matrices, Hessians, etc.) and control design methodologies (linear, nonlinear, adaptive), are calculated *on demand* by using appropriately initialized short-in time- runs of the microscopic simulator. Based on this methodology, Siettos et al. [66] present a feedback control-based framework that enables microscopic/stochastic simulators to trace their coarse-grained bifurcation diagrams past turning points by coupling the concept of coarse timesteppers with wash-out filters and the so-called "linearized" pseudo arc-length condition. In [67], Siettos et al. use an adaptive feedback control scheme to drive a microscopic simulator to its coarse critical point(s) and keep it there. The fine scale simulation is treated *as a constantly evolving experiment*; standard identification tools are wrapped around the coarse-grained observations in order to extract local, low order coarse-grained nonlinear models. Based on them, standard control design tools are used to drive the (computational) experiment to the objective.

---

[*] Corresponding author. E-mail: ksiet@mail.ntua.gr

Many other efforts aspiring to bridge microscopic simulations, experiments and control design under a macroscopic model-based perspective have been proposed. Representative examples of this interplay include the construction of stochastic partial differential equations in order to design model predictive controllers (MPC) based on kMC simulations [68, 69, 70, 71], reduction-based techniques and controller design for markovian systems [72], stochastic parameter sensitivity analysis combined with multi-step optimization [73] and the design of robust nonlinear feedforward-feedback controllers that dynamically couple kinetic Monte Carlo and finite difference simulation codes [74].

In this paper, we address a new iterative procedure based on the concept of coarse-timestepping [51, 55, 53, 56, 58, 59] that can be used to drive microscopic simulators - and potentially physical experiments- *directly* to their own open-loop unstable hyperbolic stationary points in a "vicinity" of a regular turning point. The approach does not require the implementation of a branch tracing method, such as the pseudo-arc length continuation, in order to reach the desired saddle point. If one wishes to trace the open loop codimension-one bifurcation branches of equilibria she/ he can do so in a consecutive manner, i.e. by implementing the approach for a fixed value of the bifurcation parameter and move on to the next one.

The paper is organized as follows: in section 2 we present the proposed approach for the location of unstable fixed points (in particular saddles) using the coarse-timestepper-based iterative method. In section 3 we demonstrate the scheme through two case studies: (a) a kinetic Monte Carlo simulation of a surface-reaction model and (b) an agent-based model approximating in a simplistic manner the dynamics of buying and selling of a large population of interacting individuals under the influence of mimesis. Finally, we conclude in Section 4.

## 2. Converging to a coarse-grained saddle: the procedure

### 2.1 The concept of the coarse-timestepper

Lets assume that we do not have available the explicit macroscopic equations in a closed form, but we do have an evolving microscopic (or a very-large scale) computational model. Given the (microscopic) distribution of the system $U_k \equiv U(t_k) \in R^N, N >> 1$ at time $t_k = kT$, the detailed simulator reports the values of the state variables after a time interval $T$, i.e.:

$$U_{k+1} = \wp_T(U_k, p), \qquad (1)$$

where $\wp_T : R^N \times R^m \to R^N$ is the time-evolution operator, $p \in R^m$ is the vector of systems parameters.

The main assumption behind the coarse timestepper is that a coarse-grained model for the fine-scale dynamics (1) exists and closes in terms of a few coarse-grained variables, say, $x = [x^1, x^2, ..., x^n] \in R^n, n << N$.

Usually, these are the low-order moments of the microscopically evolving distributions. The existence of a coarse-grained model implies that the higher order moments, say, $y \in R^{N-n}$, of the distribution $U$ become, relatively fast over the coarse time scales of interest, functions of the few lower ones, $x$. This scale separation can be viewed in the form of a singularly perturbed system of the following form:

$$x_{k+1} = h_s(x_k, y_k, p) \qquad (2a)$$

$$\varepsilon\, y_{k+1} = h_f(x_k, y_k, p) \tag{2b}$$

with $\varepsilon \ll 1$. Eq. (2a) corresponds to the "slow" coarse-grained dynamics while Eq. (2b) to the "fast" ones. Notice that under certain conditions [75] the overall dynamics defined by Eq. (2) will approximate very quickly the dynamics of the following system:

$$x_{k+1} = h_s(x_k, y_k, p) \tag{3a}$$

$$0 = h_f(x_k, y_k, p) \tag{3b}$$

The system above can-under some mild assumptions [75] - be written in the form of

$$x_{k+1} = F(x_k, p) \tag{4}$$

This is achieved by applying the implicit function theorem for Eq. (3b), giving:

$$y = q(x, p), \tag{5a}$$

$$F(x, p) \equiv h_s(x, q(x, p), p) \tag{5b}$$

where $q$ is a smooth continuous function defining the relation between the slow and the fast variables after a short (in the macroscopic sense) time horizon. Eq. (5a) constitutes the slow manifold on which the coarse-grained dynamics of the system evolve after a fast transient phase (see figure 1).

What a coarse timestepper does, in fact, is providing such closures *on demand*: relatively short calls of the fine scale simulation can establish this "slaving" relation (refer to [51, 55, 53, 56, 58, 59] for more detailed discussions). Briefly, once the appropriate macroscopic observables have been identified, the coarse timestepper consists of the following essential components (figure 2):

(a) Prescribe the coarse-grained initial conditions (e.g. the distribution of the temperature and concentration of the species in a combustor, the distribution of positions and velocities in a Molecular Dynamics box, e.t.c) $x(t_0) \equiv x_0$.
(b) Transform it through a lifting operator $\mu$ to consistent microscopic realizations: $U(t_0) = \mu\, x(t_0)$
(c) Evolve these realizations in time using the microscopic simulator for a short macroscopic time $T$ generating $U(t_0 + T)$. The choice of $T$ is associated with the (estimated) spectral gap of the linearization of the unavailable closed macroscopic equations [53, 60].
(d) Obtain the coarse-grained variables using a restriction operator $M$: $x_{k+1} \equiv x(t_0 + T) = MU(t_0 + T)$,

## 2.2. Description of the iterative scheme

Starting from a stable equilibrium (a node), $(x_n^*, p_i^*)$, where $x_n^* = F(x_n^*, p_i^*)$, the objective is to drive the detailed large-scale simulator to its own coarse-grained saddle point(s) (if any) *at a prescribed value of the bifurcation parameter*, say $p_i^*$, which is assumed to be "relatively close" (we explain this term below) to the turning point, say $p_i^{critical}$. We also assume that both the starting stable node and the (unknown) saddle are lying on solution branches connected through a saddle-node bifurcation at $p_i^{critical}$.

Given $p_i^*$ our procedure drives the system to the objective by creating a sequence of actions

$$u_{seq} = \{u_1, u_2, ..., u_j, ..., u_n\} \tag{6a}$$

where

$$u_j = p_i^j - p_i^* \tag{6b}$$

and $p_i^j$ is the appropriately computed value of $p_i$ at step $j$.

The sequence should converge to zero for deterministic systems (or within a small tolerance, $tol_1$, that can be defined with respect to the system's relative noise for stochastic/ microscopic systems), i.e.

$$u_n \approx 0 \tag{7a}$$

Hence, upon convergence we have that

$$p_i^j \approx p_i^* \tag{7b}$$

In a "vicinity" of the regular turning point the systems dynamics (4) close to equilibria are characterized by a separation of time scales. Hence, the trajectories $x(t)$ can be decomposed into fast $\eta(t)$ and slow $w(t)$ ones, whose dynamics can be represented as a singularly perturbed system of the form

$$\eta_{k+1} = s(\eta_k, w_k, p), \eta \in R \tag{8a}$$
$$\varepsilon w_{k+1} = f(\eta_k, w_k, p), w \in R^{n-1} \tag{8b}$$

The value of $\varepsilon$ gets smaller as the separation between the coarse-grained time-scales gets wider. After an initial time interval $t > T_r$ and under the assumptions of the Fenichel theorem [75] the dynamics of the above system can be approximated by the reduced slow system

$$\eta_{k+1} = s(\eta_k, \varphi(z_k), p) \tag{9}$$

where $\varphi$ is a sufficiently differentiable function of $\eta$, relating the fast $w$, and slow $\eta$, coarse grained variables as:

$$w = \varphi(\eta) \tag{10}$$

The above $(n-1)$ scalar equations define a curve which forms the one-dimensional slow manifold. In a neighborhood of a saddle point $(x_s, p_i)$ this is tangential to the linear unstable manifold spanned by the eigenvector $e_u$ corresponding to the dominant unstable eigenvalue $\lambda_u$, ($|\lambda_u| > 1$) of the Jacobian of the linearized system:

$$\delta x_{k+1} = \frac{\partial F}{\partial x} \delta x_k \tag{11}$$

where $\delta x_k = x_k - x_s$ and $\frac{\partial F}{\partial x}$ is the Jacobian evaluated at $x_s$.

Let us now denote the open ball in $R^n$ with radius $\varepsilon$ and center $x$ defined by

$$B(\varepsilon, x) = \{z : \|x - z\| < \varepsilon\} \tag{12}$$

In a sufficiently small neighbourhood of a saddle point $(x_s, p_i): x_s = F(x_s, p_i)$, and after some time horizon, $\forall z \in B(\varepsilon, x_s)$, excluding the points defining the coarse-grained separatrix, $\exists \delta(\varepsilon)$:

$$\|F(z, p_i) - z\| > \delta(\varepsilon) \qquad (13)$$

Eq. (13) implies that sufficiently close to the saddle point the system's dynamics will- within a certain relatively small time horizon- evolve monotonically away from the saddle, i.e.

$$\|x_{k+1} - x_s\| > \|x_k - x_s\|, \qquad (14)$$

along the (linear) one-dimensional coarse-grained unstable slow manifold. If the gap between the coarse-grained time scales in (8) is big enough, the dynamics will eventually approach very fast, and then evolve along, the "slow" eigenvector, i.e.

$$\delta x_{k+1} \to c_u \lambda_u e_u \qquad (15)$$

where the coefficient $c_u \in R$ is associated with the initial conditions. Starting from a nearby initial point at the "other side" of the separatrix, the dynamics will eventually evolve along the opposite direction, i.e. $\delta x_{k+1} \to c'_u \lambda_u (-e_u)$.

This behavior is justified by the Hartman-Grossman theorem stating that in the vicinity of hyperbolic points the local nonlinear dynamics are topologically equivalent with the linear ones [76].

In summarizing, starting from a detailed initial distribution, the system will quickly approach the "slow" coarse-grained manifold as defined by (5a) and then move along the unstable coarse-grained manifold as defined by (10); in the neighborhood of an equilibrium point and under the above assumptions, the coarse-grained dynamics will evolve along $e_u$ i.e., along the eigenvector corresponding to the dominant eigenvalue of the linearized problem (11) (see figure 3).

Depending on the initial conditions $x_0$ across the coarse-grained separatrix, the system will ultimately approach the stable fixed point $x_n$ or move towards another attractor, say $A$ (such as another stable fixed point).

Due to continuity, for small perturbations of the control parameter, i.e. for $u_{j+1} = u_j + \delta u$, the direction of the dominant eigenvector as well as the value of the corresponding eigenvalue and generally the behaviour of the trajectories do not change considerably , i.e. they remain qualitatively unchanged. Hence by monitoring an observable-such as the difference $(x^i_{k+1} - x^i_k)$ with the maximum norm- we can infer about the side of the separatrix that the coarse grained dynamics evolve.

Under the above assumptions our algorithm reads as follows:

**Step 0**. Initialize the sequence with $u_0 > 0$ or $u_0 < 0$ if $p_i^* < p_i^{critical}$ and $p_i^* > p_i^{critical}$ respectively. A loose estimate of the value of $p_i^{critical}$ can be obtained by plain temporal simulations.

As we consider the parameters' values of the vector $\{p - \{p_i\}\}$ as constants, we rewrite the coarse timestepper model on the slow manifold as

$$x_{k+1} = F(x_k, u) \qquad (16)$$

**Step I.** Initialize the algorithm with the required parameters: $gain, tol_1, tol_2, c, tol_3$.

All parameters of the iterative protocol may be in principle constants or can be changed adaptively. Set $k = 1$ and go to step II.

**Step II.** ***Do While convergence to the saddle point*** *(i.e. when* $|u_k| < tol_1$*)*

**A.** We "perturb" the system for a macroscopic time horizon $T$ using $u_{k-1}$ as our "control" action to get

$$U_k = \wp_T(U_{k-1}, u_{k-1}) \tag{17}$$

**B.** Set the initial interval of the possible values of $u_k$ as follows:

$$\begin{aligned}
&\text{If } u_{k-1} > 0 \text{ Then} \\
&\qquad \text{set interval} = [lB_0, uB_0] = [0, u_{k-1}]; \\
&\text{Else} \\
&\qquad \text{set interval} = [lB_0, uB_0] = [u_{k-1}, 0]; \\
&\text{Endif}
\end{aligned}$$

Set $q = 0$;

**C. *Do while*** $|h(d)| > tol_2$ *(see below, how* $h(d)$ *is defined) using a Chord-like technique*

1. *Update* $q = q + 1$; compute the next "control" action from

$$u_k = lB_{q-1} + a(uB_{q-1} - lB_{q-1}) \tag{18}$$

where $0 < a < 1$.

2. Compute the coarse-grained state vector by letting the detailed simulator to relax on the coarse-grained slow manifold corresponding to the new value of $u_k$. This can be achieved by running the detailed model for a sufficiently small time horizon $T_r$:

$$U_{k+1} = \wp_{T_r}(U_k, u_k) \tag{19}$$

Upon the coarse-grained slow manifold, monitor the evolution of the dynamics for a short time $T_s$:

$$x_{T_s} = F_{T_s}(x_0, u_k), \tag{20}$$

where $x_0$ is the coarse-grained state vector as obtained by restricting $U_{k+1}$.

3. Calculate the function $h : R^n \to R$ reading

$$h(x_{T_s}) = x^i_{T_s} - x^i_0, \ |x^i_{T_s} - x^i_0| := max^n_{j=1}(|x^j_{T_s} - x^j_0|) \tag{21}$$

4. Update the search space of $u_k$ according to the following manner:

    a. *If* $h(d) < 0$ (meaning -under our convention- that the coarse-grained fixed point of (20) is located at $p_i > p_i^* + u_k$) Then set $lB_q = u_k$;

    *Else If* $h(d) > 0$ (meaning that the equilibrium is located at $p_i < p_i^* + u_k$) Then set $uB_q = u_k$;

*Endif*

b. *If* ($uB_0 == uB_q$ or $lB_0 == lB_q$ ) and $uB_q - lB_q < tol_3$; *i.e. the procedure has minimized the search space to some extend and the control action cannot drive the system from one region of attraction to the other Then*

    *If* $uB_0 == uB_q$; *Then update the bounds of the control action as:*

$$uB_0 = uB_0 + c, \ uB = uB_0;$$

*Else*

$$lB_0 = lB_0 - c, \ lB = lB_0;$$

*Endif*
*Endif*

### *End Do While (of Step II. C)*

**D.** Update $k = k + 1$ and run the detailed model for a time horizon $T_r$ : $\boldsymbol{U}_k = \wp_{T_r}(\boldsymbol{U}_{k-1}, u_{k-1})$ in order to relax the system to the slow manifold.

**E.** *Update the new action as*

$$u_k = gain \cdot u_{k-1}, \tag{22}$$

The value of the gain should be greater than one in order to drive the system (applying the control action as defined by (18)) closer to the seeking saddle point (see figure 4). Its value could be constant or adaptive depending on how far is the current control action $p_i^j$ from the nominal value $p_i^*$. If the system is currently at the coarse-grained saddle point, $u_{k-1} = 0$, the next action will be also $u_k = 0$.
Set $k = k + 1$.

### *End Do While (of Step II)*

    The procedure finds the intermediate "control/iterative" actions through the chord method in order to "traverse" in some way the turning point and drive the system to the saddle point that we are seeking.

    A feature of the method is that it does not require the reconstruction -through lifting- of consistent high-dimensional microscopic/detailed distributions from the low-dimensional coarse-observables. What is fed into the algorithm is a "snapshot" of the microscopic distribution ($U_k$ in Eq. (19) ) whose projection $\boldsymbol{x}_0$ into the low-dimensional space serves as the initial guess into the chord-like iteration [77].

    Under our assumptions the steering of the system towards the coarse-grained saddle is achieved by manipulating the "control" parameter by-passing the necessity for using a conventional Newton-Chord-like technique [77]. Hence, it does not require the numerical approximation of the system's coarse-grained Jacobian $\nabla \boldsymbol{F}(\boldsymbol{x}_k)$, which would (a) increase the computational cost requiring the computation of all first-order partial derivatives of the vector field $\boldsymbol{F}(\cdot)$, (b) require the derivation of consistent to the coarse-grained variables detailed/ microscopic distributions.

    For stochastic systems Step II.C.2 should be repeated for a number of ensembles and then take the average of the results.

## 3. Case Studies and Simulation Results

We demonstrate our method using two case studies. The first one is a kMC realization of a chemical reaction on a catalytic surface while the second one is a simplistic agent-based model of buying and selling dynamics under the influence of mimesis which bears strong analogies to integrate-and-fire models of neurons.

### 3.1 Kinetic Monte Carlo simulations of surface reactions

We first present an example describing the dynamics of simple catalytic reactions using stochastic simulations. The species react, are adsorbed or desorbed on a finite lattice with periodic boundary conditions. At each time instant, the sites of the lattice are considered to be either vacant or occupied by the reaction species. In the kMC context the system's dynamics are described by the following chemical master equation:

$$\frac{dP(x,t)}{dt} = \sum_y Q(y,x)P(y,t) - \sum_y Q(y,x)P(x,t) \tag{23}$$

where $P(x,t)$ is the probability that the system will be in state $x$ at time $t$. $Q(y,x)$ denotes the probability for transition from state $y$ to $x$ per unit time. The summation runs over all possible transitions (reactions), say, $N_r$.

Here we used the Gillespie kMC algorithm [78, 79] to numerically simulate the above stochastic equation. The basic steps read as follows:

*Do While { time < $t_{end}$ }*

    Step 1. Calculate the cumulative transition probability $Q_0(N_r) = \sum_{y=1}^{N_r} Q(y,x)$.

    Step 2. Draw two random numbers, say $r_1$ and $r_2$, from a uniform distribution in the interval [0 1].

    Step 3. Select the reaction $y^*$ from the set of all possible events so that the following inequality is satisfied:

$$Q_0(y^* - 1) < r_1 Q_0(N_r) < Q_0(y^*) \tag{24}$$

    Step 4. Perform the selected event and update all kinetic rates according to the occurred event.

    Step 5. Advance time by $\Delta t = -\frac{1}{Q_0(N_r)} \ln r_2$ , time = time + $\Delta T$

*End DoWhile*

Our illustrative microscopic model is a kMC realization of a simplification of the dynamics of catalytic CO oxidation:

$$CO + \frac{1}{2}O_2 \rightarrow CO_2 \tag{25}$$

The reaction mechanism can be schematically described by the following elementary steps:

$$\begin{aligned}&(1)\ CO_{gas} + *_i \leftrightarrow CO_{ads,i} \\ &(2)\ O_{2,gas} + *_i + *_j \leftrightarrow O_{ads,i} + O_{ads,j} \\ &(3)\ CO_{ads,i} + O_{ads,j} \rightarrow CO_{2,gas} + *_i + *_j\end{aligned} \tag{26}$$

where $i$, $j$ are sites on the square lattice, * denotes a site with a vacant adsorption site, while "ads" denotes adsorbed particles. A simple schematic of the above mechanism is given in figure 5.

The macroscopic mean field description is given by the following equations [61]:

$$\dot{\theta}_A = \alpha(1 - \theta_A - \theta_B) - \gamma\theta_A - 4k_r\theta_A\theta_B \qquad (27a)$$

$$\dot{\theta}_B = 2\beta(1 - \theta_A - \theta_B)^2 - 4k_r\theta_A\theta_B \qquad (27b)$$

where $\theta_l$ represent the coverage of species ($l$ =A, B, corresponding to $CO$ and $O$) on the catalytic surface. The parameters $a$, $\beta$, $\gamma$ are associated with $CO$ adsorption, $O$ dissociative adsorption and $CO$ desorption rates while $k_r$ is a reaction rate constant.

Here, the system's size and the number of realizations were chosen to be $N_{size} = 1000 x 1000$ and $N_{runs} = 2000$, respectively. The value of the time horizon was selected as $T = 0.05$.

The coarse timestepper

$$\theta_{A,\ k+1} = \Phi_{1,T}(\theta_A, \theta_B, \beta, a, k_r) \qquad (28a)$$
$$\theta_{B,\ k+1} = \Phi_{2,T}(\theta_A, \theta_B, \beta, a, k_r) \qquad (28b)$$

of the kMC realization was used as a "black box" timestepper.

Figure 6a depicts the convergence of the coarse-grained variables for $a = 1.6$, $\beta = 4$, $\gamma = 0.04$ $k_r = 1$, to the correct value of the saddle point at $\theta_A^* = 0.7323$, $\theta_B^* = 0.0881$ starting from either ($\theta_A^* = 0.9701$, $\theta_B^* = 0.0016$), or ($\theta_A^* = 0.1126$, $\theta_B^* = 0.6902$) stable equilibria; figure 6b depicts the corresponding "control" action sequence of the parameter $a$.

By repeatedly calls of the procedure for different values of the parameter $\beta$ we also calculated the corresponding coarse-grained bifurcation diagram with respect to $\beta$ (shown in figure 7). The obtained diagram coincides with the one obtained through the deterministic system (22) (see also [61]).

To this end, figure 8 depicts the diagram of the two eigenvalues of the Jacobian of the linearized deterministic system (22) at the fixed points with respect to the bifurcation parameter $\beta$, justifying the method's assumption on the time-scale separation in the coarse-grained dynamics.

### 3.2 Agent-based simulation of a financial market caricature

Our second illustrative example is an agent-based model describing a simplistic buyer-seller procedure appearing in [80]. The state $x_i$ of each individual $i = 1, 2, ..., N$ is a scalar in the interval $(-1\ 1)$. The value of $x_i$ reflects the "mood" of agent $i$ to buy or sell. At -1 individuals sell, at +1 they buy; then their mood is reset to neutral (zero). Thus, the main activity is the evolution of the state $x_i$ of the individuals within $(-1\ 1)$; this is modelled by the following stochastic differential equation:

$$\frac{dx_i(t)}{dt} = -\gamma x_i(t) + I_i(t) \qquad (29)$$

coupled with the constraint that when the state reaches the thresholds -1, 1 it is immediately reset to $x_i = 0$. $I_i(t)$ is a stochastic variable representing the incoming information arriving to the individual $i$ at time $t$. This can be considered as the exogenous information each individual draws from its environment (e.g. mass media news, or opinions of stock market consultants) plus the influence due to recent decisions of other agents. Each agent reacts to such "good" or "bad" news, and its state, $x_i$, is increased by $\varepsilon^+$ or decreased by $-\varepsilon^-$, respectively. The arrival of the news comes at discrete times with Poisson distributed rates of $v_{ex}^+$ and $v_{ex}^-$ respectively. Hence, the state of each agent changes in an event-driven way by a series of positive or negative jumps at times $t^{(k)}, k = 1, 2, 3, \ldots$. In the absence of such stimuli, the state of each individual is exponentially attracted to zero (no particular interest to buy or sell) at rate $\gamma$. Agents are also influenced by the overall buying and selling rates. These rates are $R^+$ and $R^-$, defined as the fraction of agents which buy or sell per unit time. The apparent reaction rates to good and bad news are

$$v^\pm = v_{ex}^\pm + gR^\pm \tag{30}$$

respectively, where $g$ is a factor representing the relative influence of the buying and selling rates. These rates are defined as the average rates over a small, finite interval, which serves as our reporting system time horizon, say $\delta t$.

For our simulations, we used the following values for the parameters: $v_{ex}^+ = 20$, $v_{ex}^- = 20$, $v_{ex}^+ = 20$, $v_{ex}^- = 20$, $\varepsilon^+ = 0.075$, $\varepsilon^- = 0.075$, $\gamma = 1$, $\delta t = 0.25$. The number of agents of the system was set to $N = 50000$; the observable of our method is the average state $\bar{x} = \frac{1}{N}\sum_{i=1}^{N} x_i$.

Figure 9 depicts the coarse-grained bifurcation diagrams of $\bar{x}$ vs. $v_{ex}^\pm$ as constructed by repetitive calls of the proposed procedure for $g = 38$, while figure 10 shows the coarse-grained bifurcation diagram $\bar{x}$ with respect to the parameter $g$ as this obtained by implementing the protocol in an iterative way for consecutive values of $g$. Figures 11a,b show the convergence of $\bar{x}$ and the computed "control-iterative" actions required to drive the system at the upper coarse-grained saddle stationary point for $g = 38$, using $v_{ex}^+$ as the "control" parameter. Figures 11c, d illustrate the convergence of $\bar{x}$ to the lower saddle point and the "control" actions sequence, respectively, using $v_{ex}^-$ as the "control" parameter.

## 4. Conclusions

Microscopic simulation is a vital tool in the study of complex and multi-scale phenomena ranging from Material Science and Mechanics to Engineering, Ecology, Biology and Social Sciences.

We introduced and demonstrated a new iterative procedure based on the concept of coarse-timestepping [53] that can be used to drive microscopic simulators to their own unknown open-loop coarse-grained saddle points at a chosen value of the bifurcation parameter when no macroscopic/ coarse-grained equations are available.

This is accomplished through the construction of a sequence of "control" actions, using the bifurcation parameter as an "actuator", and designing its dynamics utilizing basic numerical bifurcation theory concepts. The procedure builds the actions through a chord-like method in order to drive the microscopic simulator around its own coarse-grained criticality and finally steer it to the coarse-grained saddle point that is sought for a given configuration of the parameter space.

The procedure is illustrated using two microscopic simulations, namely (a) a kinetic Monte Carlo implementation of a simple surface reaction scheme for which the macroscopic equations for the expected values of the state variables were known so that we could validate our results, and, (b) an agent-based model describing in a simplistic way the dynamics of many interacting "investors" under mimesis. Both models exhibit coarse regular turning points marking the onset of coarse-grained saddles.

Using the microscopic simulators as black box coarse-timesteppers we were able to locate the coarse-grained saddles at specified values of the bifurcation parameters. We also constructed in an iterative way the expected open-loop coarse-grained bifurcation diagrams.

The methodology described here is not only applicable to multi-scale/ microscopic systems but also to very-large systems as those arise in many fields such as material science, chemical and combustion engineering and biological problems. It may possibly serve as a basis for the development of new adaptive control protocols that can be used to drive physical experiments directly to their own coarse-grained saddles.

The basic prerequisites of the approach are (a) the occurrence of a time-scale separation so that the detailed dynamics can be decomposed into fast and coarse grained slow ones, (and if this holds) (b) the a-priori knowledge of the coarse-grained "slow" variables, and finally (c) a clear separation between time-scales between the unstable and stable coarse-grained modes (which is valid in a neighbourhood of regular turning points).

The first assumption implies that a coarse-grained model for the dynamics at the macroscopic/continuum level exists and closes in terms of a few coarse-grained variables. However, in several complex problems these are not known before-hand. In this case, one should first resort to advanced techniques for data reduction such as Principal Component Analysis, or Diffusion Maps [81] that can suggest the right coarse-grained observables. Further research could also proceed towards the connection of the method with techniques such as the Computational Singular Perturbation (CSP) [82] for the systematic numerical reconstruction of the slow manifolds. In that way, one would systematically converge to the slow manifold in the coarse-grained phase space and thus relax the third assumption of the method.

Acknowledgements: This work was partially supported by the State Scholarships Foundation of Greece (IKY) and the National Technical University of Athens through the Basic Research Program ''Constantine Caratheodory''.

**Figures**

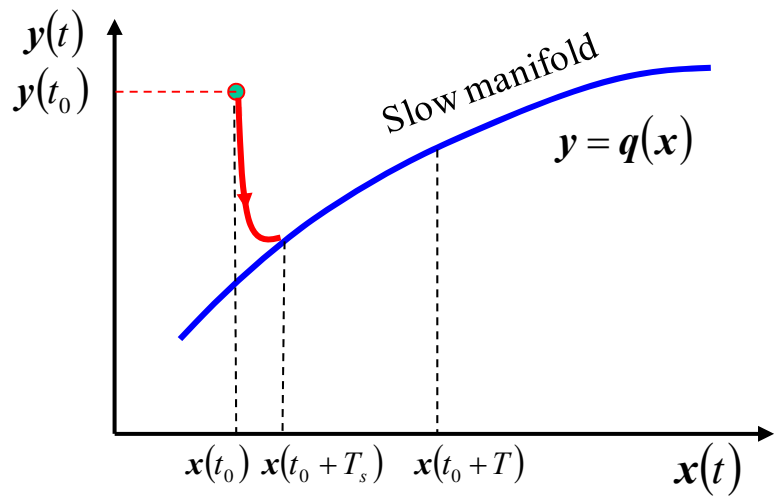

**Figure 1.**

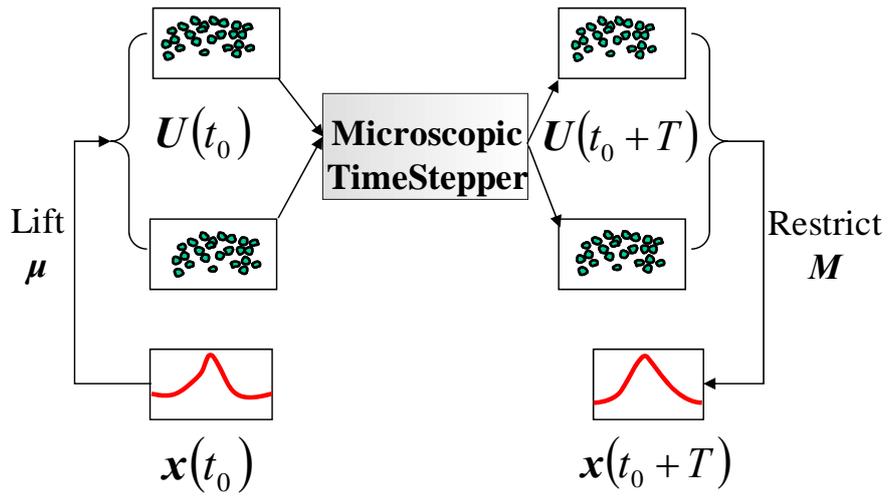

**Figure 2.**

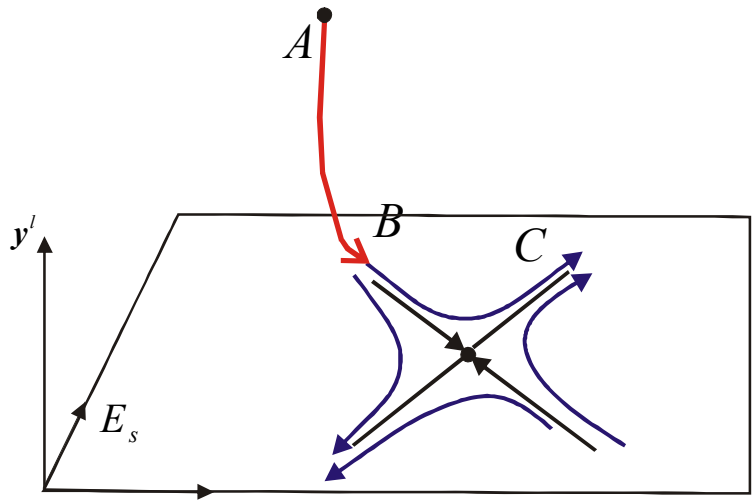

**Figure 3.**

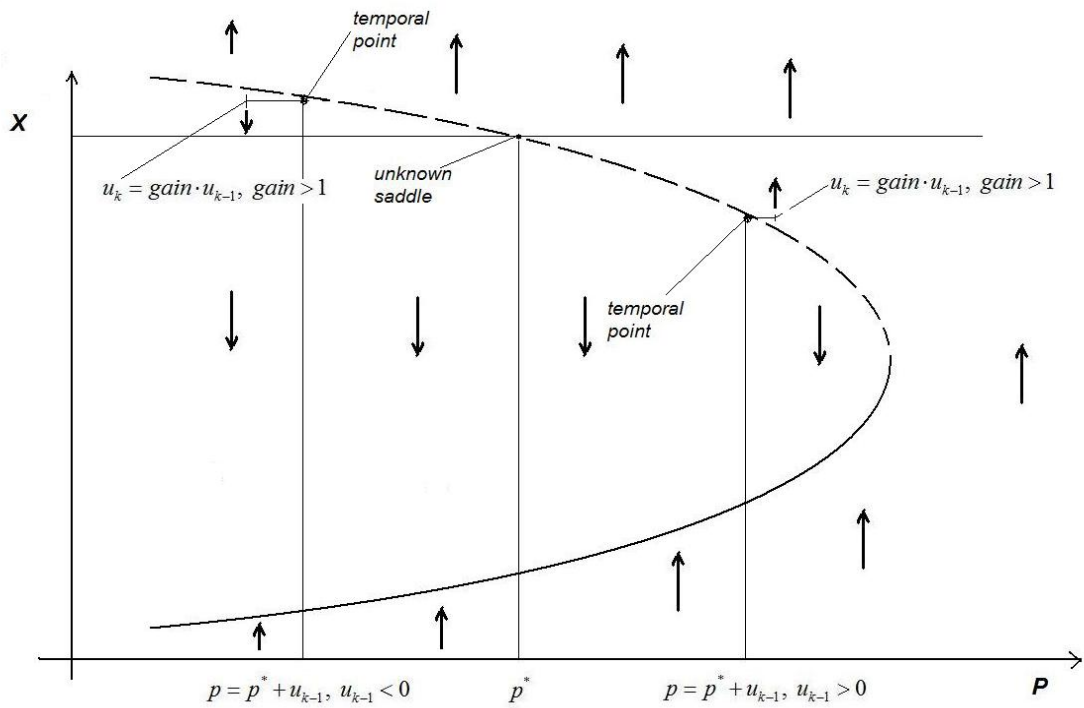

**Figure 4.**

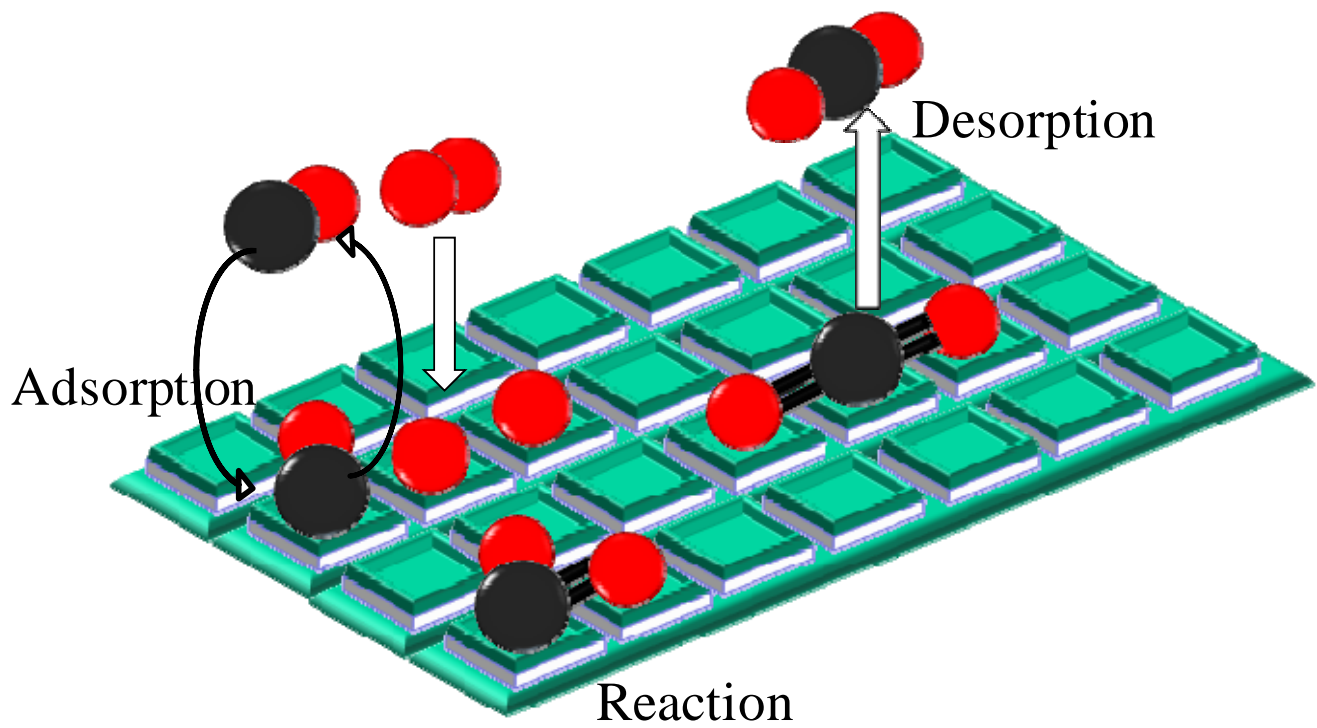

**Figure 5.**

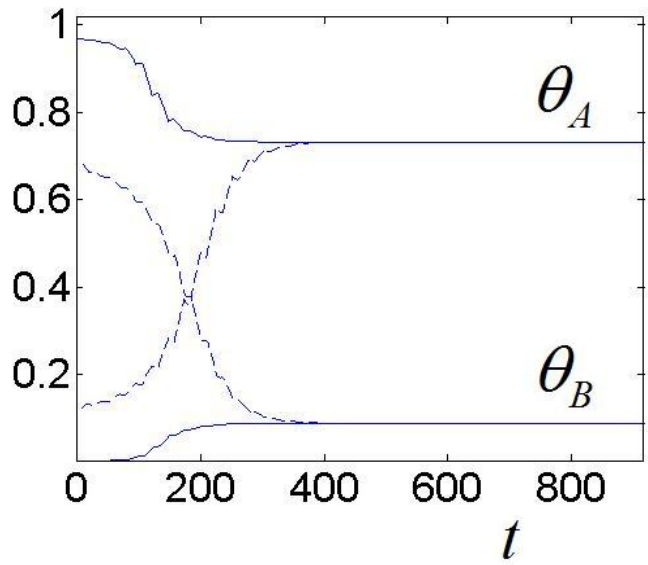

(a)

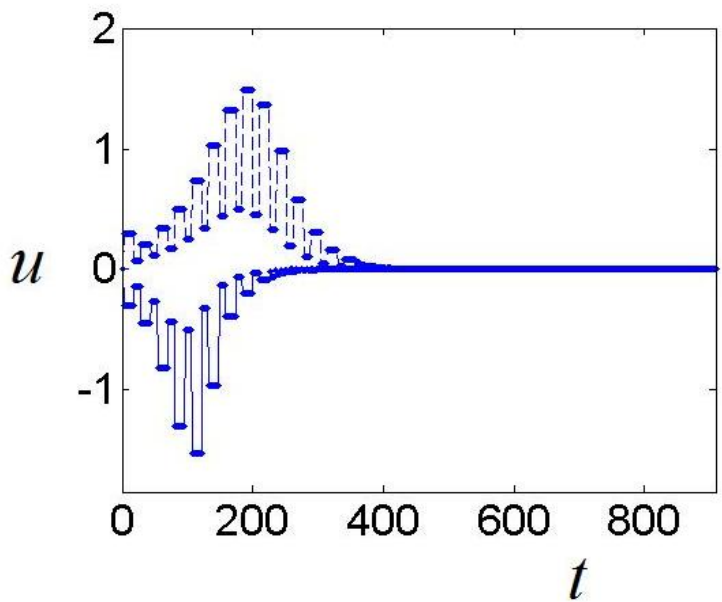

(b)

**Figure 6.**

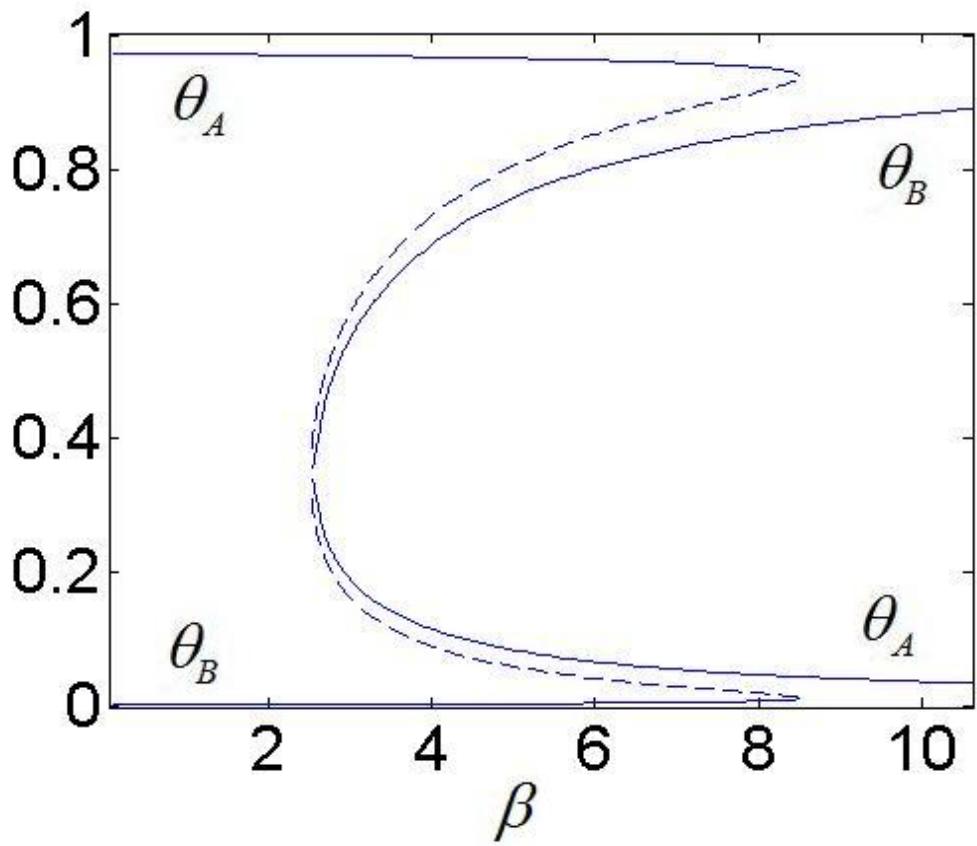

**Figure 7.**

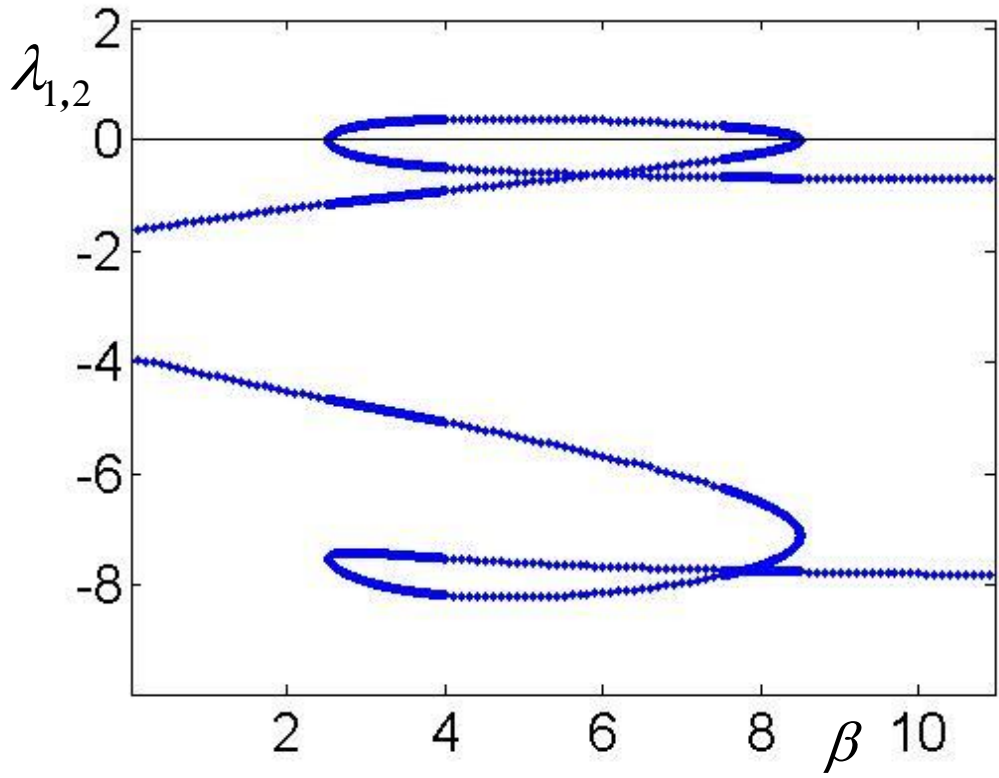

**Figure 8.**

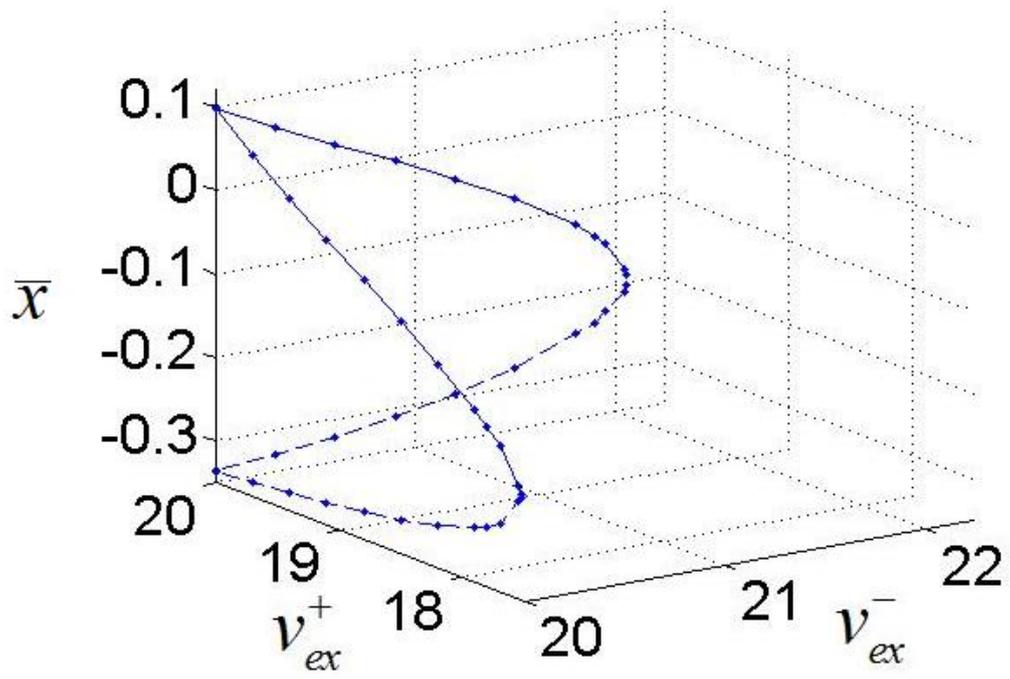

**Figure 9.**

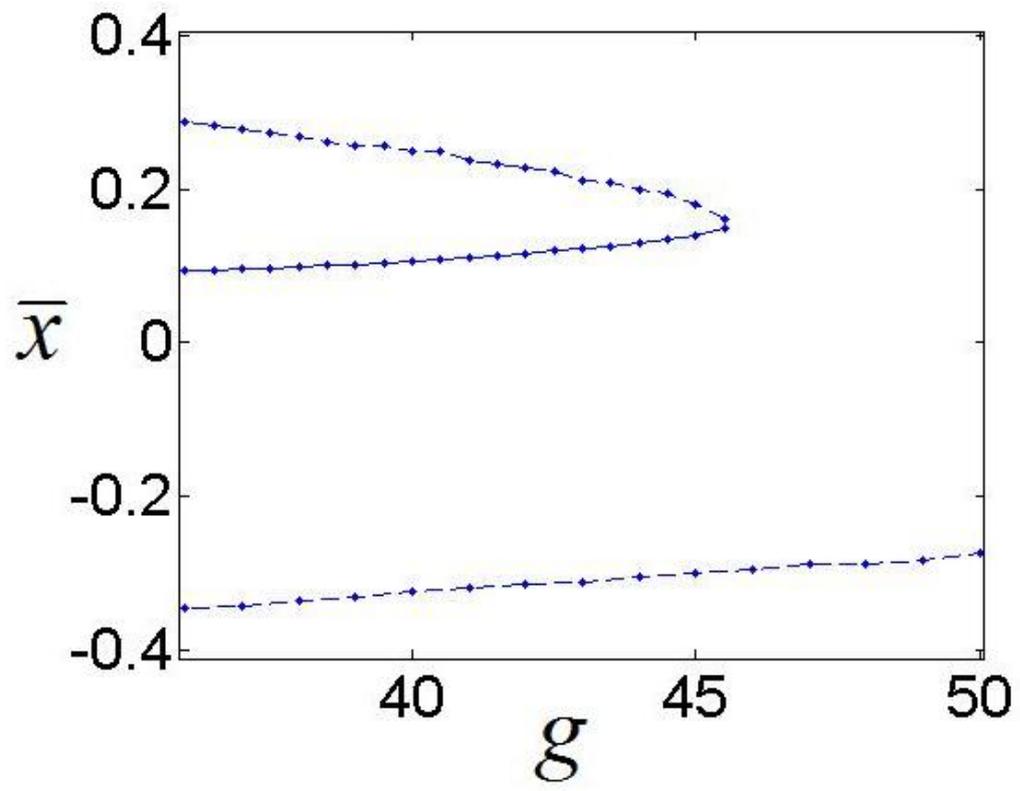

**Figure 10**

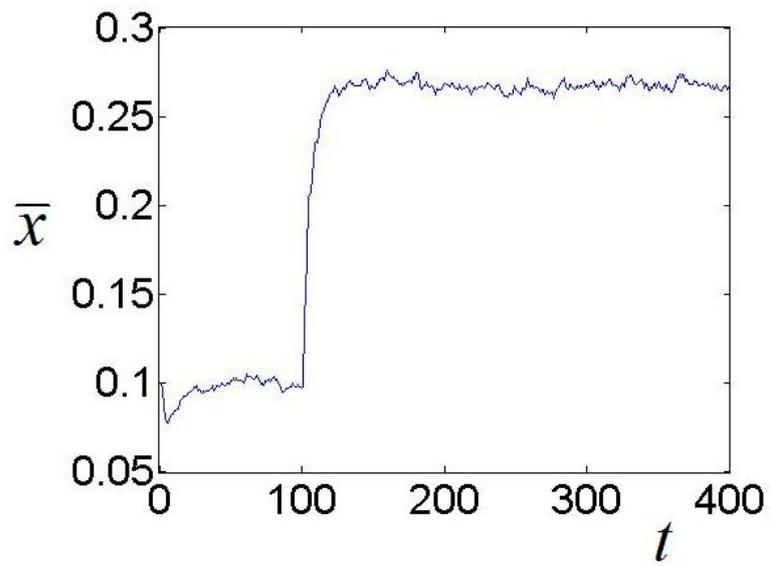

(a)

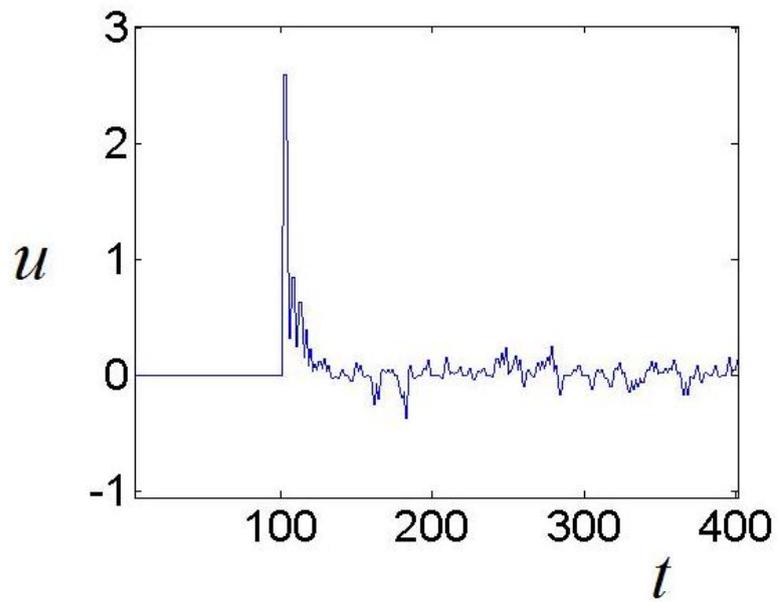

(b)

**Figure 11**

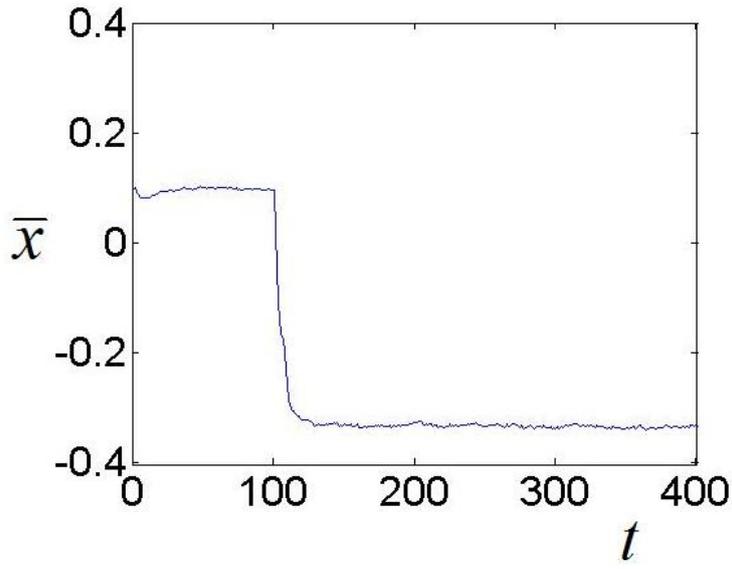

(c)

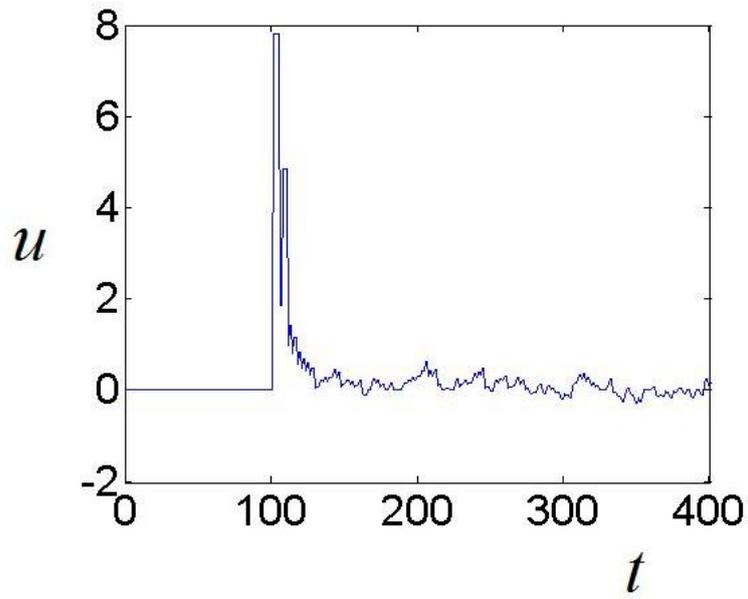

(d)

**Figure 11**

**Captions of Figures**

**Figure 1.** The basic assumption of the Equation-Free approach: The higher order moments $y$ of an evolving microscopic distribution become "quickly"- within the time-horizon $T$ - functionals of the lower-order moments $x$.

**Figure 2.** A schematic of the coarse-timestepper.

**Figure 3.** Trajectories in the neighborhood of a coarse-grained saddle point. Starting from a detailed distribution (point A) the system will, under the assumptions, relax on the slow coarse-grained manifold, $E_s$ (AB orbit). Then the coarse-grained dynamics will evolve towards the coarse-grained unstable manifold (BC orbit).

**Figure 4**. A saddle-node bifurcation diagram. The use of *gain*>1 places the current point into the "right" basin of attraction in order to move the system closer to the seeking saddle point.

**Figure 5.** A simple schematic of the CO oxidation on a catalytic surface.

**Figure 6. (a)** The trajectories and (b) the "control" actions sequences using $\alpha$ as "control" parameter to the saddle point for $a = 1.6$, $\beta = 4$, $\gamma = 0.04$ $k_r = 1$.

**Figure 7**. Coarse-grained bifurcation diagram of the catalytic CO oxidation with respect to the parameter $\beta$.

**Figure 8**. The two eigenvalues corresponding to the Jacobian of the deterministic mean field model (27) at the fixed points. The "slow" eigenvalue changes sign at the turning points. There is a clear separation in time scales.

**Figure 9**. The coarse-grained bifurcation diagrams of $\bar{x}$ vs. $v_{ex}^-$ and $v_{ex}^+$ as constructed by repetitive calls of the proposed approach at $g = 38$ for different values of $v_{ex}^-$ and $v_{ex}^+$ respectively.

**Figure 10**. The coarse-grained bifurcation diagram of the agent-based model, with $g$ serving as the bifurcation parameter. It was obtained by implementing the proposed iterative protocol for consequent values of $g$.

**Figure 11a-d.** The trajectories and "control" actions sequences to the upper saddle point, (figures a-b), and to the lower saddle point, (figures c-d), starting from the stable node at $g = 38$ using $v_{ex}^+$ and $v_{ex}^-$ respectively.